\documentclass[11pt]{article}
\usepackage{amsmath}
\usepackage{amsmath,amssymb,latexsym,color}
\usepackage[mathscr]{eucal}
\newtheorem{thm}{Theorem}[section]

\newtheorem{lemma}[thm]{Lemma}

\newtheorem{example}{Example}[section]
\newtheorem{defin}{Definition}[section]

\newtheorem{remark}{Remark}[section]

\newcommand{\proof}{{\it Proof.\quad}}
\newcommand{\qed}{\hfill\Box\medskip}

\usepackage{CJK}
\begin{document}
\begin{CJK*}{GBK}{song}

\newcommand{\be}{\begin{equation}\label}
\newcommand{\ee}{\end{equation}}
\newcommand{\bea}{\begin{eqnarray}\label}
\newcommand{\eea}{\end{eqnarray}}

\title{The chromatic spectrum of 3-uniform bi-hypergraphs
}

\author{
Ping Zhao$^{\rm a}$\quad Kefeng Diao$^{\rm a}$\quad    Kaishun
Wang$^{\rm b}$\thanks{Corresponding
author: wangks@bnu.edu.cn}\\
{\footnotesize a. \em School of Science, Linyi University,
Linyi, Shandong, 276005, China }\\
{\footnotesize b. \em  Sch. Math. Sci. {\rm \&} Lab. Math. Com.
Sys., Beijing Normal University, Beijing 100875,  China} }

\date{}
 \maketitle

\begin{abstract}

   Let $S=\{n_1,n_2,\ldots,n_t\}$ be a finite set
 of positive integers with
$\min(S)\geq 3$ and $t\geq 2$. For any positive integers $s_1,s_2,\ldots,s_t$,
 we construct a family of 3-uniform  bi-hypergraphs ${\cal H}$
 with the feasible set $S$ and $r_{n_i}=s_i, i=1,2,\ldots,t$, where each $r_{n_i}$ is
 the $n_i$th component of the chromatic spectrum of ${\cal H}$. As a
 result, we solve one open problem for $3$-uniform bi-hypergraphs proposed
   by Bujt\'{a}s and Tuza   in 2008.
Moreover,   we find a family of sub-hypergraphs  with the same
feasible set and the same chromatic spectrum as it's own. In
particular,  we obtain a small upper bound on the minimum number of
vertices in 3-uniform bi-hypergraphs with any given feasible set.

\medskip
\noindent {\em AMS Subject classification:} 15A36

\noindent {\em Key words:} Bi-hypergraph; feasible set; chromatic
spectrum
\end{abstract}

\section{Introduction}

A {\em mixed hypergraph } on a finite set $X$ is a triple ${\cal
H}=(X, {\cal C}, {\cal D})$, where ${\cal C}$ and ${\cal D}$ are
families of subsets of $X$, called the {\em ${\cal C}$-edges} and
{\em ${\cal D}$-edges}, respectively. A set $B\in {\cal C}\cap {\cal
D}$ is called a \emph{bi-edge} and a \emph{bi-hypergraph} is a mixed
hypergraph with ${\cal C}={\cal D}$, in which case we may simply
write ${\cal H}=(X, {\cal B})$ instead of  ${\cal H}=(X, {\cal B,
\cal B})$. In general, we follow the terminology of
\cite{Voloshin2}.

A sub-hypergraph ${\cal H}'=(X', {\cal C}', {\cal D}')$ of a mixed
hypergraph ${\cal H}=(X, {\cal C}, {\cal D})$  is a {\em partial
sub-hypergraph}  if $X'=X$, and
 ${\cal H}'$ is called a {\em derived sub-hypergraph}
  of  ${\cal H}$ on $X'$, denoted by ${\cal H}[X']$, when ${\cal C}'=\{C\in {\cal C}| C\subseteq X'\}$
 and ${\cal D}'=\{D\in {\cal D}| D\subseteq X'\}$.  If $|C|=|D|=r\geq 2$  for any $C\in {\cal C}$ and $D\in {\cal D}$,
  then the mixed hypergraph ${\cal H}=(X, {\cal C}, {\cal D})$ is {\em $r$-uniform}.

  Two mixed hypergraphs ${\cal H}_1=(X_1, {\cal C}_1, {\cal D}_1)$
and ${\cal H}_2=(X_2, {\cal C}_2, {\cal D}_2)$ are \it isomorphic
\rm if there exists a bijection $\phi$ between $X_1$ and $X_2$ that maps
each $C$-edge of ${\cal C}_1$ onto a $C$-edge of ${\cal C}_2$ and
maps each $D$-edge of ${\cal D}_1$ onto
 a $D$-edge of ${\cal D}_2$,
 and vice versa. The bijection $\phi$ is called an \emph{isomorphism} from ${\cal H}_1$ to ${\cal H}_2$.

  A {\em proper $k$-coloring} of
${\cal H}$ is a mapping from $X$ into a set of $k$ {\em colors} so
that each ${\cal C}$-edge has two vertices with a {\em Common} color
and each ${\cal D}$-edge has two vertices with {\em Distinct}
colors. A {\em strict $k$-coloring} is a proper $k$-coloring using
all of the $k$ colors, and a mixed hypergraph is {\em $k$-colorable}
if it has a strict $k$-coloring. The maximum (minimum) number of
colors in a strict coloring of ${\cal H}=(X, {\cal C}, {\cal D})$ is
the {\em upper chromatic number} $\overline \chi ({\cal H})$
 (resp. {\em lower chromatic number} $\chi ({\cal H})$) of ${\cal H}$.
The study of   the colorings of mixed hypergraphs has made a lot of
progress since its inception \cite{Voloshin}.
 For more information, see \cite{Tuza}.

  The set of all the  values $k$ such that ${\cal H}$ has a strict $k$-coloring is called the {\em feasible set}
  of ${\cal H}$, denoted by $\Phi({\cal H})$. For each $k$, let $r_k$ denote
the number of {\em partitions} of the vertex set. Such partitions
are called {\em feasible partitions}.  The vector $R({\cal
H})=(r_1,r_2,\ldots,r_{\overline\chi})$ is called the {\em chromatic
spectrum} of ${\cal H}$.

It is readily seen that if $1\in \Phi({\cal H})$, then ${\cal H}$ cannot have any ${\cal D}$-edges.
 For the case $1\notin \Phi({\cal H})$, Jiang  et al.     \cite{Jiang}  proved
 that,
 for any finite set $S$ of integers greater than 1, there exists a mixed hypergraph ${\cal H}$
such that $\Phi({\cal H})=S$, and Kr$\acute{\rm a}$l \cite{Kral}
strengthened this result by showing that prescribing any positive
integer   $r_k$, there exists a mixed hypergraph which has precisely
$r_k$ $k$-coloring for all $k\in S$. Recently, Bujt\'{a}s and Tuza
\cite {Bujtas1}
    gave the necessary and sufficient condition for a finite set $S$
    of natural numbers being the feasible set of an $r$-uniform mixed
    hypergraph, and they raised the following two open problems:

\medskip  {\bf Problem 1.}  Determine the chromatic
      spectrum of  $r$-uniform bi-hypergraphs.

\medskip   {\bf Problem 2.} Determine the minimum number of
vertices in
    $r$-uniform bi-hypergraphs with   given feasible set.

\medskip
 The motivation of this paper is to solve these open problems. We focus on 3-uniform bi-hypergraphs and organize this paper as follows. In Section 2, for any integer $t\geq 2$,
 any finite set  $S=\{n_1,n_2,\ldots,n_t\}$
 of integers with
$\min(S)\geq 3$ and any positive integers $s_1,s_2,\ldots,s_t$, we
construct a family of 3-uniform bi-hypergraphs ${\mathcal H}$  with
the feasible set $S$ and $r_{n_i}=s_i, i\in [t]$, where $r_{n_i}$ is
the $n_i$th component of the chromatic spectrum $R({\cal H})$. As a
result, we solve Problem 1 for $3$-uniform bi-hypergraphs. Moreover,
we discuss the maximality of the bi-hypergraph $\mathcal H$ in terms
of number of bi-edges, i.e., by adding any other bi-edge its
chromatic spectrum changes. In Section 3, we find a family of
sub-hypergraphs of ${\mathcal H}$ with the same feasible set and the
same chromatic spectrum as it's own. In particular, we give a small
upper bound on the minimum
 number of vertices in $3$-uniform bi-hypergraphs with given feasible set.

\section{Construction}

For any positive integer $n$, let $[n]$ denote the set
$\{1,2,\ldots,n\}$.
We first introduce the construction.\\

For any  integers $s\geq 2$ and $n_1\geq
\cdots \geq n_s\geq 3$, let $X_{n_1,\ldots,
n_s}=\{(x_1,\ldots,x_s)|x_j\in [n_j], j\in [s]\}$ and
$${\cal
B}_{n_1,\ldots,
n_s}=\{\{(x_1,\ldots,x_s),(y_1,\ldots,y_s),(z_1,\ldots,z_s)\}|~|\{x_j,
y_j, z_j\}|=2, j\in [s]\}.$$ Then $(X_{n_1,\ldots, n_s},{\cal B}_{n_1,\ldots ,n_s})$
is a 3-uniform bi-hypergerpah, denoted by ${\cal H}_{n_1,\ldots,n_s}$.

 Note that, for any $i\in [s]$, $c_i^s=\{X_{i1}^s,\ldots,X_{in_i}^s\}$ is a strict
$n_i$-coloring of ${\cal H}_{n_1,\ldots,n_s}$,
  where $$X_{ij}^s=\{(x_1,x_2,\ldots,x_{i-1},j, x_{i+1},\ldots,x_s)|x_k\in [n_k],
k\in [s]\setminus \{i\}\}, j\in [n_i].$$ In the following we will
prove   that $c_1^s,\ldots,c_s^s$ are all the strict colorings of
${\cal H}_{n_1,\ldots,n_s}$ by induction on $s$.

\begin{lemma} For any integers  $n_1> n_2\geq 3$, we have
$$\Phi({\cal H}_{n_1,n_2})=\{n_1,n_2\} \mbox{~and~} r_{n_1}=r_{n_2}=1.$$
 \end{lemma}

\proof
  Suppose $c=\{C_1,C_2,\ldots, C_m\}$ is a strict coloring of ${\cal H}_{n_1,n_2}$.
  Then there are the following  two possible cases:

  \medskip

 \textbf{Case 1} $(1,1)$ and $(1,2)$ fall into a common color class.

  Suppose $(1,1),(1,2)\in C_1$.
The bi-edge $\{(1,1),(1,2),(2,1)\}$ implies $(2,1)\notin C_1$.
Suppose  $(2,1)\in C_2$. Since $\{(1,1),(2,1),(2,2)\}$ and $\{(1,1),(1,2),(2,2)\}$ are bi-edges,
 $(2,2)\in C_2$. Similarly, we have $(i,1),(i,2),\ldots, (i,n_2)\in C_i$ for each
 $i\in [n_1]$. It follows that
    $c=c_1^2$.

      \medskip

 \textbf{Case 2} $(1,1)$ and $(1,2)$ fall into distinct color classes.

 Similar to Case 1,  we have $c=c_2^2$.

 Hence, the desired result follows.
 $\qed$

\begin{thm} For any integers $s\geq 2$ and $n_1> n_2> \cdots >n_s\geq
3$, we have   $$\Phi({\cal H}_{n_1,\ldots,n_s})=\{n_1,\ldots, n_s\}
\mbox{~and~} r_{n_1}=\cdots =r_{n_s}=1.$$
\end{thm}

\proof By Lemma 2.1,  the conclusion is true for $s=2$.

Let $X'=\{(x_1,x_3,x_3,x_4,\ldots,x_s)|x_j\in [n_j],
j\in [s]\setminus \{2\}\}$.
Then ${\cal H}'={\cal H}_{n_1,\ldots,n_s}[X']$ is isomorphic to ${\cal
H}_{n_1,n_3,n_4,\ldots,n_s}$. By induction, all the strict colorings of ${\cal H}'$
are as follows: $$c'_i=\{X'_{i1},X'_{i2},\ldots,
X'_{in_i}\}, i\in [s]\setminus \{2\},$$ where
$X'_{ij}=X'\cap X_{ij}^s, j\in [n_i]$.

For any strict coloring $c=\{C_1,C_2,\ldots, C_m\}$
 of ${\cal H}_{n_1,\ldots,n_s}$, there are the
following three possible cases:

   \medskip

\textbf{Case 1}  $c|_{X'}=c_1'$.

For each $i\in [n_1]$ and $x_j\in [n_j], j\in [s]\setminus \{1\}$, the bi-edges
\begin{eqnarray*}&&\{(i,x_2,x_3,\ldots,x_s),(i,x_3',x_3',x_4',\ldots,x_s'),(i',x_3',x_3',x_4',\ldots,x_s')\},\\
  &&\{(i,x_2,x_3,\ldots,x_s),(i,x_3',x_3',x_4',\ldots,x_s'),(i'',x_3',x_3',x_4',\ldots,x_s')\}
\end{eqnarray*}
 imply that
 $(i,x_2,x_3,\ldots,x_s)\in C_i$.  Hence, $c=c_1^s$.

  \medskip

\textbf{Case 2} $c|_{X'}=c_3'$.

From the bi-edge $\{(1,1,2,1,\ldots,1),(1,\ldots,1),(2,\ldots,2)\}$,  we get
 $(1,1,2,1,\ldots,1)\in C_1\mbox{~or~}C_2$.

    \medskip

\textbf{Case 2.1 }
$(1,1,2,1,\ldots,1)\in C_1$.

For any $j\in [s]\setminus \{2\}$ and $x_j\in [n_j]$, the bi-edges
\begin{eqnarray*}&&\{(x_1,1,x_3,x_4,\ldots,x_s),(x_1',x_3',x_3',x_4',\ldots,x_s'),(1,1,2,1,\ldots,1)\},\\
  &&\{(x_1,1,x_3,x_4,\ldots,x_s),(x_1',x_3'',x_3'',x_4',\ldots,x_s'),(1,1,2,1,\ldots,1)\}
\end{eqnarray*}
imply that $(x_1,1,x_3,x_4,\ldots,x_s)\in C_1$,  where $|\{x_j,x_j',1\}|=2$
for any $j\in [s]\setminus \{2,3\}$ and $|\{x_3, x_3', 2\}|=|\{x_3, x_3'', 2\}|=2, x_3''\neq x_3'$.
For any $j\in [s]\setminus \{2\}, x_j\in [n_j]$ and $k\in [n_3]\setminus \{1\}$,
from the bi-edges
\begin{eqnarray*}&&\{(x_1,k,x_3,x_4,\ldots,x_s),
(x_1,1,x_3,x_4,\ldots,x_s),(x_1',1,x_3',x_4',\ldots,x_s')\},\\ &&\{(x_1,k,x_3,x_4,\ldots,x_s),
(x_1,k,k,x_4',\ldots,x_s'),(x_1',1,x_3',x_4,\ldots,x_s)\},
\end{eqnarray*} we have  $(x_1,k,x_3,x_4,\ldots,x_s)\in C_k$, where $|\{k, x_3, x_3'\}|=2$.

 For any $k\in [n_3]$, the bi-edge $$\{(x_1,n_3+1,x_3,x_4,\ldots,x_s),(x_1,k,x_3,x_4,\ldots,x_s),
 (x_1',k,x_3',x_4',\ldots,x_s')\}$$ implies that $(x_1,n_3+1,x_3,x_4,\ldots,x_s)\notin C_{k}$.
   Suppose  $(1,n_3+1,1,\ldots,1)\in C_{n_3+1}$. Then
  the bi-edge $\{(1,n_3+1,1,\ldots,1),(x_1,n_3+1,x_3,x_4,\ldots,x_s),(x_1',1,x_3',x_4',\ldots,x_s')\}$ implies that
   $(x_1,n_3+1,x_3,x_4,\ldots,x_s)\in C_{n_3+1}$, where
   $|\{1,x_i,x_i'\}|=2, x_i, x_i'\in [n_i],i\in [s]\setminus \{2\}$.
   Similarly, for
 any $j\in [n_2-n_3]$, we have $(x_1,n_3+j,x_3,x_4,\ldots,x_s)\in C_{n_3+j}$.
 Therefore,  $c=c_2^s$.

   \medskip

\textbf{Case 2.2 } $(1,1,2,1,\ldots,1)\in C_2$.

 Similar to Case 2.1, we have
    $c=c_3^s$.

       \medskip

\textbf{Case 3} There exists a $k\in [s]\setminus \{1,2,3\}$ such that $c|_{X'}=c_k'$.

Similar to Case 1, we obtain $c=c_k^s$.

Hence, the desired result follows. $\qed$

By Theorem 2.2, we get the the following  result:

\begin{thm} For any integer $t\geq 2$, any finite set $S=\{n_1, n_2, \ldots, n_t\}$ of integers
with $\min (S)\geq 3$ and any positive integers $s_1, s_2,
\ldots, s_t$, ${\cal
H}_{n_1,\ldots,n_1,n_2,\ldots,n_2,\ldots, n_t,\ldots,n_t}$  is a 3-uniform bi-hypergraph  with
feasible set
 $S$ and $r_{n_i}=s_i, i\in [t]$,  where each $n_i$  appears $s_i$ times in the vector
  $(n_1,\ldots,n_1,n_2,\ldots,n_2,\ldots, n_t,\ldots,n_t)$.
\end{thm}

This theorem answers Problem 1 proposed by Bujt\'{a}s and Tuza
\cite{Bujtas1}.

The following theorem discover the maximality of the bi-hypergraph
${\cal H}_{n_1,\ldots,n_s}$ in terms of number of bi-edges.

\begin{thm} For any integers $t\geq 2$, $n_1>n_2>\cdots >n_t\geq 3$ and  $s_1,s_2,\ldots,s_t$, let $${\cal
H}=(X_{n_1,\ldots,n_1,\ldots,n_t,\ldots,n_t}, {\cal
B}_{n_1,\ldots,n_1,\ldots,n_t,\ldots,n_t}\cup \{B\}),$$ where  $B\notin {\cal
B}_{n_1,\ldots,n_1,\ldots,n_t,\ldots, n_t}.$ Then
$$R({\cal
H})\neq R({\cal H}_{n_1,\ldots,n_1,\ldots,n_t,\ldots,n_t}).$$
\end{thm}

\proof Let $s=s_1+\cdots + s_t$. Note that ${\cal
H}_{n_1,\ldots,n_1,\ldots,n_t,\ldots,n_t}$ is a partial sub-hypergraph of ${\cal H}$,
 so ${\cal H}$ has no other strict colorings except $c_1^s,\ldots,c_s^s$.

Suppose $B=\{(a_1,\ldots,a_s), (b_1,\ldots,b_s),
(c_1,\ldots,c_s)\}$. If there exists an $m\in [s]$ such that
$|\{a_m,b_m,c_m\}|=3$, then $c_m^s$ is not a strict coloring of
${\cal H}$ since it makes the three vertices of $B$ to be colored
with three distinct colors. Otherwise, there exists an $m\in [s]$
such that $|\{a_m,b_m,c_m\}|=1$. It follows that $c_m^s$  makes the
three vertices of $B$ to be colored with only one color. Hence, the
desired result follows.$\qed$

\section{Sub-hypergraphs}

In this section,  we  find a family of sub-hypergraphs of ${\cal
H}_{n_1,\ldots,n_s}$  with the same feasible set and the same chromatic spectrum as its own.

For any integers $s\geq 2$ and $n_1\geq n_2>\ldots>n_s>3$,
let
\begin{eqnarray*}X_1&=&\{(n_2+j,1,\ldots,1),(n_2+j,n_2,n_3,\ldots,n_s)|j\in [
n_1-n_2]\}\\
X_i&=&\bigcup_{j=1}^{n_i-n_{i+1}}\{(n_{i+1}+j,1,\ldots,1),(n_{i+1}+j,\ldots,n_{i+1}+j,n_{i+1},\ldots,n_s)\}\\
&\cup&(\bigcup_{j=0}^{n_i-n_{i+1}}\{(1,n_{i+1}+j,\ldots,n_{i+1}+j,1,\ldots,1)\}), i=2,\ldots,s-1,\\
X_s&=&(\bigcup_{i=1}^3\bigcup_{k=1}^3\{(i,k,\ldots,k)\})
\cup (\bigcup_{k=4}^{n_s}\{(1,k,\ldots,k),(k,1,\ldots,1),(k,\ldots,k)\})
 \end{eqnarray*} and
$$X_{n_1,\ldots, n_s}^*=\bigcup_{i=1}^sX_i,~~~~
{\cal H}_{n_1,\ldots,n_s}^*={\cal
H}_{n_1,\ldots,n_s}[X_{n_1,\ldots,n_s}^*],$$ where $n_{i+1}+j$
appears $i-1$ times in $(1,n_{i+1}+j,\ldots,n_{i+1}+j,1,\ldots,1)$.
The vertex $(1,n_{i+1},\ldots,n_{i+1},1,\ldots,1)$ is called the
\emph{inflexion} of $n_{i+1}$.

In the rest we  shall prove that ${\cal H}_{n_1,\ldots,n_s}^*$ is a
family of the desired sub-hypergraphs of ${\cal
H}_{n_1,\ldots,n_s}$.

For any
 $i\in [s]$, $c_i^{s*}=\{X_{i1}^{s*},\ldots,X_{in_i}^{s*}\}$ is a strict
$n_i$-coloring of ${\cal H}^*_{n_1,\ldots,n_s}$,
  where
$X_{ij}^{s*}=X_{n_1,\ldots, n_s}^*\cap X_{ij}^s, j\in [n_i].$

\begin{lemma} For any integers $n_1\geq n_2>3$, we have
$$\Phi({\cal H}^*_{n_1,n_2})=\Phi({\cal H}_{n_1,n_2}) \mbox{~and~}R({\cal H}^*_{n_1,n_2})=R({\cal H}_{n_1,n_2}).$$
\end{lemma}

\proof Suppose $c=\{C_1,C_2,\ldots,C_m\}$ is a strict coloring of
${\cal H}^*_{n_1,n_2}$. We  get the following two possible cases:

 \medskip

\textbf{Case 1}  $(i,1),(i,2),(i,3)\in C_i, i=1,2,3$.

 For any
$j\in [n_2]\setminus \{1,2,3\}$, the bi-edges $\{(1,j),
(1,1),(2,1)\}$ and $\{(1,j),(1,1),(3,1)\}$ imply that $(1,j)\in
C_1$. From the bi-edges $\{(4,1),(1,1),(1,2)\}, \{(4,1),(2,1),
(2,2)\}$ and $\{(4,1),(3,1),(3,2)\}$, we have $(4,1)\notin C_1\cup
C_2\cup C_3$. Suppose $(4,1)\in C_4$. Then since  $\{(4,1),(1,1),
(4,4)\}$ and $\{(1,1),(1,4),(4,4)\}$ are bi-edges, we have $(4,4)\in
C_4$. Similarly, $(i,1),(i,i)\in C_i$ for each $i\in [n_2]\setminus
\{1,2,3\}$. If $n_1> n_2$,  then $(n_2+j,1),(n_2+j,n_2)\in C_{n_2+j}$ for
each  $j\in [n_1-n_2]$. Hence,  $c=c_1^{2*}$.

  \medskip

\textbf{Case 2} $(1,j),(2,j),(3,j)\in C_j, j=1,2,3.$

Similar to Case 1, we have $c=c_2^{2*}$.

Hence, the desired result follows. $\qed$

\begin{thm} For any integers $s\geq 2$ and $n_1\geq n_2>\cdots > n_s> 3$, we have
$$\Phi({\cal H}^*_{n_1,\ldots,n_s})=\Phi({\cal H}_{n_1,\ldots,n_s}),
R({\cal H}^*_{n_1,\ldots,n_s})=R({\cal H}_{n_1,\ldots,n_s}).$$
\end{thm}

\proof By Lemma 3.1,  the conclusion is true for $s=2$.

Let $X^*=X_{n_1,\ldots, n_s}^*\setminus (X_2\cup X_1)$. Then ${\cal H}^*={\cal H}_{n_1,\ldots,n_s}[X^*]$
is isomorphic to ${\cal H}^*_{n_3,n_3,n_4,\ldots,n_s}$.
 By induction,  all the strict colorings of ${\cal H}^*$ are as follows:
 $$c_i^*=\{X_{i1}^*,X_{i2}^*,\ldots,X_{in_i}^*\}, i\in [s]\setminus \{2\},$$
 where $X_{ij}^*=X_{ij}^{s*}\cap X^*, j\in [n_i]$.

Suppose $c=\{C_1,C_2,\ldots,C_m\}$ is a strict coloring of ${\cal
H}^*_{n_1,\ldots,n_s}$. We  have the following three possible cases:

  \medskip

\textbf{Case 1} $c|_{X^*}=c_1^*$.

 For any $j\in [n_2-n_3]\cup \{0\}$, since
\begin{eqnarray*}&&\{(1,2,\ldots,2),(2,2,\ldots,2),(1,n_3+j,1,\ldots,1)\},\\ &&\{(1,3,\ldots,3),
(3,3,\ldots,3),(1,n_3+j,1,\ldots,1)\}
\end{eqnarray*} are bi-edges,
$(1,n_3+j,1,\ldots,1)\in C_1$. Further, for each $j\in [n_1-n_3]$,
the bi-edges
\begin{eqnarray*}&&\{(n_3+j,1,\ldots,1),(1,1,\ldots,1),(1,2,\ldots,2)\},\\
&&\{(n_3+j,1,\ldots,1),(k,1,\ldots,1),(k,k,\ldots,k)\} \mbox{~or~}\\
&&\{(n_3+j,1,\ldots,1),(k,1,\ldots,1),(k,\ldots,k,n_p,\ldots,n_s)\}
 \end{eqnarray*} imply that  $(n_3+j,1,\ldots,1)\notin C_k$ for any $k\in [n_3]$,
  where $k\in [n_s]$ or $n_p<k\leq n_{p-1}$ for some $p\in [s]\setminus \{1,2,3\}$. Suppose $(n_3+1,1,\ldots, 1)\in C_{n_3+1}$.
  From the bi-edges
 \begin{eqnarray*}&&\{(n_3+1,1,\ldots, 1), (n_3+1,n_3+1,n_3,\ldots,n_s), (k,1,\ldots,1)\}, \\
 &&\{(n_3+1,1,\ldots, 1), (n_3+1,n_3+1,n_3,\ldots,n_s), (k',1,\ldots,1)\},
 \end{eqnarray*} we have  $(n_3+1,n_3+1,n_3,\ldots,n_s)\in
C_{n_3+1}$. Similarly, $(n_3+j,1,\ldots,1),
(n_3+j,n_3+j,n_3,\ldots,n_s)\in C_{n_3+j}$ for each $j\in
[n_2-n_3]\setminus \{1\}$. Further, if $n_1> n_2$, then
$(n_2+j,1,\ldots,1), (n_2+j,n_2,\ldots,n_s)\in C_{n_2+j}$ for each
$j\in [n_1-n_2]$. Therefore,  $c=c_1^{s*}$.

  \medskip

\textbf{Case 2} $c|_{X^*}=c_3^*.$

Note that  $(n_3+j,1,\ldots,1)\in C_1$ for any $j\in [n_1-n_3]$. Further,
the bi-edge $\{(1,n_3,1,\ldots,1),
(n_3,n_3,n_3,n_4,\ldots,n_s),(1,\ldots,1)\}$ implies that
$(1,n_3,1,\ldots,1)\in C_1$ or $C_{n_3}$.

 \medskip

\textbf{Case 2.1}
$(1,n_3,1,\ldots,1)\in C_1$.

For any $j\in [n_2-n_3]$, from  the bi-edges
\begin{eqnarray*}&&\{(1,n_3,1,\ldots,1),
(n_3,n_3,n_3,n_4,\ldots,n_s),(1,n_3+j,1,\ldots,1)\},\\
&&\{(1,n_3,n_3,1,\ldots,1), (n_3,n_3,n_3,n_4,\ldots,n_s),
(1,n_3+j,1,\ldots,1)\}
\end{eqnarray*} and
\begin{eqnarray*}&&\{(1,n_3,n_3,1,\ldots,1), (n_3+j,n_3+j,n_3,\ldots,n_s),
(1, n_3, 1,\ldots,1)\},\\
&&\{(1,1,1,\ldots,1),(n_3+j,n_3+j,n_3,\ldots,n_s),
(1,n_3+j,1,\ldots,1)\},
\end{eqnarray*} we have $(1,n_3+j,1,\ldots,1)\in C_1$ and $(n_3+j,n_3+j,n_3,\ldots,n_s)\in
C_{n_3}$. Similarly, if $n_1> n_2$, then for any $j\in [n_1-n_2]$,
$(n_2+j,n_2,n_3,\ldots,n_s)\in C_{n_3}$. Hence,  $c=c_3^{s*}$.

 \medskip

\textbf{Case 2.2}
$(1,n_3,1,\ldots,1)\in C_{n_3}$.

 For any $k\in [n_3]\setminus \{1\}$, the bi-edge
 \begin{eqnarray*}&&\{(1,k,\ldots,k),
 (k,k,\ldots,k),(1,n_3+1,1,\ldots,1)\} \mbox{~or~}\\
  &&\{(1,k,\ldots,k,1,\ldots,1),(k,\ldots,k,n_p,\ldots,n_s),(1,n_3+1,1,\ldots,1)\}
  \end{eqnarray*}
  imply that $(1,n_3+1,1,\ldots,1)\notin C_k$, where $k\in [n_s]$ or $n_p<k\leq n_{p-1}$ for some $p\in [s]\setminus \{1,2,3\}$. Further, from the bi-edges
 \begin{eqnarray*}&&\{(n_3+1,n_3+1,n_3,\ldots,n_s),
 (1,n_3+1,1,\ldots,1),(1,n_3,n_3,1,\ldots,1)\},\\ &&\{(n_3+1,n_3+1,n_3,\ldots,n_s),
 (1,n_3+1,1,\ldots,1),(n_3+1,1,\ldots,1)\},\\
 &&\{(n_3+1,n_3+1,n_3,\ldots,n_s), (1,n_3,1,\ldots,1),(1,n_3,n_3,1,\ldots,1)\},
  \end{eqnarray*} we have $(1, n_3+1,1,\ldots,1)\notin C_1$. Suppose $(1,n_3+1,1,\ldots,1)\in C_{n_3+1}$. The bi-edges
  \begin{eqnarray*}&&\{(n_3+1,n_3+1,n_3,\ldots,n_s),(1,n_3+1,1,\ldots,1),
  (1,n_3,1,\ldots,1)\},\\  &&\{(n_3+1,n_3+1,n_3,\ldots,n_s),(1,n_3,1,\ldots,1),(1,n_3,n_3,1,\ldots,1)\}
   \end{eqnarray*} imply that $(n_3+1,n_3+1,n_3,\ldots,n_s)\in C_{n_3+1}$.
   Similarly, for each $j\in [n_2-n_3]\setminus \{1\}$, $(1,n_3+j,1,\ldots,1),
   (n_3+j,n_3+j,n_3,\ldots,n_s)\in C_{n_3+j}$. Further, if $n_1>n_2$, we have
   $(n_2+j,n_2,n_3,\ldots,n_s)\in C_{n_2}$ for each $j\in [n_1-n_2]$. Hence,
$c=c_2^{s*}$.

 \medskip

\textbf{Case 3} There exists a $k\in [s]\setminus \{1,2,3\}$ such
that $c|_{X^*}=c_k^*$.

Note that $(n_3+j,1,\ldots,1)\in C_1$ for any  $j\in [n_1-n_3]$. For any $ j\in [n_2-n_3]\cup \{0\}$, the
bi-edges\begin{eqnarray*}&&\{(n_k,\ldots,n_k,n_{k+1},\ldots,n_s), (1,n_k,\ldots,n_k,1,\ldots,1),(1,n_3+j,1,\ldots,1)\},\\
   &&\{(n_k,\ldots,n_k,n_{k+1},\ldots,n_s), (1,n_k,\ldots,n_k,1,1,\ldots,1),(1,n_3+j,1,\ldots,1)\}
  \end{eqnarray*} imply that  $(1,n_3+j,1,\ldots,1)\in C_1$. From the bi-edges
   \begin{eqnarray*}&&\{(1,n_3+j,1,\ldots,1),(n_3+j,1,\ldots,1),(n_3+j,n_3+j,n_3,\ldots,n_s)\},\\
   &&\{(1,n_k,\ldots,n_k,1,\ldots,1), (1,n_k,\ldots,n_k,1,1,\ldots,1),(n_3+j,n_3+j,n_3,\ldots,n_s)\},
   \end{eqnarray*} we have $(n_3+j,n_3+j,n_3,\ldots,n_s)\in C_{n_k}$ for any $ j\in [n_2-n_3]$,
    where the vertex $(1,n_k,\ldots,n_k,1,1,\ldots,1)$ is the inflexion of $n_k$. Similarly, if $n_1> n_2$, then
     $(n_2+j,n_2,n_3,\ldots,n_s)\in C_{n_k}$ for any $j\in [n_1-n_2]$. Therefore,  $c=c_k^{s*}$.

Hence, the desired result follows. $\qed$

 Note that $|X^*_{n_1,\ldots,n_s}|=2n_1+n_2+s-2$. Hence
 $2n_1+n_2+s-2$ is the upper bound of   the minimum number of vertices in 3-uniform bi-hypergraphs
with any feasible set $\{n_1,n_2,\ldots, n_s\}$.

\section*{Acknowledgment}

The research
 is supported by NSF of
Shandong Province (No. ZR2009AM013), NCET-08-0052, NSF of China
(10871027) and the Fundamental Research Funds for the Central Universities of China.

\end{CJK*}


\begin{thebibliography}{00}
\bibitem{Bujtas1}C. Bujt$\acute{\rm a}$s and Zs. Tuza, Uniform mixed hypergraphs: the possible numbers of colors,  Graphs Combin. 24 (2008) 1-12.

%\bibitem{Bujtas2}C. Bujt$\acute{\rm a}$s and Zs. Tuza, Color-bounded hypergraphs, I: General results, Discrete Math. 309 (2009) 4890-4902.

\bibitem{Kral}D. Kr$\acute{\rm a}$l, On feasible sets of mixed hypergraphs,    Electronic
J. Combin.  11 (2004)  $\sharp$R19.

\bibitem{Jiang}T. Jiang, D. Mubayi, Zs. Tuza, V. Voloshin and D.
West, The chromatic spectrum of mixed hypergraphs,  Graphs Combin.
18 (2002) 309-318.

\bibitem{Tuza} Zs. Tuza and V. Voloshin, Problems and results on
colorings of mixed hypergraphs, Horizons of Combinatorics, Bolyai
Society Mathematical Studies 17, Springer-Verlag, 2008, 235-255.

\bibitem{Voloshin}V. Voloshin, On the upper chromatic number of a
hypergraph,   Australasian J. Combin.   11 (1995) 25-45.

\bibitem{Voloshin2} V. Voloshin,  Coloring Mixed Hypergraphs: Theory, Algorithms
and Applications, AMS, Providence, 2002.


 \end{thebibliography}
\end{document}